\documentclass[11pt]{amsart}

\usepackage{amsmath,amssymb,amsthm}

\usepackage{fancyhdr}
\usepackage{amsfonts,graphicx}

\theoremstyle{plain}
\newtheorem{Theo}{Theorem}[section]
\newtheorem{lem}[Theo]{Lemma}
\newtheorem{prop}[Theo]{Proposition}

\theoremstyle{plain}
\theoremstyle{definition}

\theoremstyle{remark}

\newtheorem*{rema*}{Remark}

\newcommand{\ZZ}{\mathbb{Z}}  
  
\newcommand{\NN}{\mathbb{N}}

\newcommand{\RR}{\mathbb{R}}

\newcommand{\QG}{(\textnormal{QG})_{\frac{1}{2}}}
\newcommand{\QGG}{(\textnormal{QG})_{\alpha}}

\date{}

\title[the critical 2D  quasi-geostrophic equation]
{On the global well-posedness of the critical  quasi-geostrophic equation}
\author[H. Abidi]{Hammadi Abidi}
\address{IRMAR, Universit\'e de Rennes 1\\ Campus de
Beaulieu\\ 35~042 Rennes cedex\\ France}
\email{hamadi.abidi@univ-renne1.fr}
\author[T. Hmidi]{Taoufik Hmidi}
\address{IRMAR, Universit\'e de Rennes 1\\ Campus de
Beaulieu\\ 35~042 Rennes cedex\\ France}
\email{thmidi@univ-renne1.fr}

\begin{document}
\maketitle

\begin{abstract}
We prove the global well-posedness  of the critical  dissipative
quasi-geostrophic equation for large initial data  belonging to the critical Besov space $\dot B^0_{\infty,1}(\RR^2).$
\end{abstract}

\section{Introduction}
In this paper we are concerned  with the initial value problem of the $2$D dissipative quasi-geostrophic equation
$$
\QGG\left\lbrace
\begin{array}{l}\partial_{t}\theta+v\cdot\nabla\theta+|\hbox{D}|^\alpha\theta=0\\
\theta_{\vert t =0}=\theta^0,
\end{array}
\right.$$
where  the scalar function $\theta$ represents the potential temperature and the parameter $\alpha\in]0,2]$. The velocity 
$v=(v^1,v^2)$ is determined by Riesz transforms of $\theta:$
$$
v=(-\partial_{2}|\hbox{D}|^{-1}\theta,\partial_{1}|\hbox{D}|^{-1}\theta)
:=(-{R}_{2}\theta,{R}_{1}\theta);\quad|\hbox{D}|=\sqrt{-\Delta}.
$$
In addition to its intrinsic mathematical importance this equation serves as  a $2$D model 
in geophysical fluid dynamics, for more details about the subject \mbox{see \cite{C-M-T,Ped}}.   

This equation   has been intensively investigated and  much attention is carried to the problem of global well-posedness.  For the sub-critical case
$(\alpha>1)$ the theory seems to be in a satisfactory state. Indeed, global existence and uniqueness for arbitrary  initial data  are established in various function spaces \mbox{(see for example  
\cite{C-W,Res}).} However the  critical  and super-critical cases, corresponding respectively to $\alpha=1$ and $\alpha<1,$ are harder to deal with. In the super-critical case, we have until now only global results for small initial data, see for instance \cite{Cha,Chen,Keraani,Ning1,Wu1,Wu2}. For critical case , Constantin, C\'ordoba and Wu  showed in \cite{C-C-W} the global existence in Sobolev space $H^1$   under smallness  assumption of  $L^\infty$ norm of $\theta^0.$ Many other relevant results can be found in \cite{C-C,Ju,Ju1,marchand}. Very recently, Kiselev, Nazarov and Volberg  proved in \cite{volberg} the global well-posedness  for arbitrary periodic smooth initial data by using an elegant argument of modulus of continuity. In \cite{Caffarelli}, Caffarelli and Vasseur established the global regularity of  weak solutions associated to $L^2$ initial data.

The main goal of this work is to establish the global well-posedness   in the critical case   when initial data belong to  the homogeneous critical Besov space $\dot{B}_{\infty,1}^0 (\RR^2):$ we remove the periodic condition and we weaken the initial regularity. Before giving our main result let us first precise the notion of critical spaces.   Let $\theta$ be  a solution of $\QG$ and $\lambda>0$ 
 then $\theta_{\lambda}(t,x)=\theta(\lambda t,\lambda x)$ is also a solution. One class of  scaling invariant spaces is the homogeneous Besov spaces  $(\dot B_{p,r}^{\frac{2}{p}}),$ \mbox{with 
 $p,\,r\in[1,\infty].$}
    
Our first main result reads as follows. 
\begin{Theo}\label{Thm1}
Let $\theta^0\in\dot B^0_{\infty,1},$ then there exists a unique global solution $\theta$ to $\QG$ such that
$$
\theta\in {\mathcal{C}}(\RR_+;\,\dot B^0_{\infty,1})
\cap L^1_{\textnormal loc}(\RR_+;\,\dot B^1_{\infty,1}).
$$
\end{Theo}

The proof relies essentially on two facts: the first one is the establishement of the local existence which is the major part of this paper and the derivation of some smoothing effects of the solution, described in next Theorem. The second one is the use of the modulus of continuity as done in \cite{volberg}. We mention that the property allowing us to remove the periodicity is the spatial decay of the solution. The key of our local existence result is some new estimates for the following transport-diffusion equation~:
$$
(\textnormal{TD})\left\lbrace
\begin{array}{l}
\partial_t \theta+v\cdot\nabla \theta+|\textnormal{D}|\theta=f\\
{\theta}_{| t=0}=\theta^{0},\\
\end{array}
\right.
$$
where $v$ and $f$ are given and $\theta$ is the unknown  scalar function . We state now  our second main result.

\begin{Theo}\label{Thm3}
Let $s\in]-1,1[,\, r,\overline r\in[1,+\infty]$ with $r\geq \bar{r},$ $f\in\widetilde L^{\overline r}_{\textnormal{loc}}(\RR_{+}; \dot{B}_{\infty,1}^{s+\frac{1}{\overline r}-1})$ and $v$ be a divergence free  vector field belonging to $L^1_{\textnormal{loc}}(\RR_{+};\textnormal{Lip}(\RR^2)).$ 
We consider a smooth solution $\theta$ of the transport-diffusion equation $(\textnormal{TD}),$
then there exists a constant $C$ depending only on $s$   such that for every  $t\in\RR_{+}$
 $$
 \|\theta\|_{\widetilde L^r_{t}\dot{B}_{\infty,1}^{s+\frac{1}{r}}}\leq C e^{C\int_{0}^t\|\nabla v(\tau)\|_{L^\infty}d\tau}
 \big(  \|\theta^0\|_{\dot{B}_{\infty,1}^s}+\|f\|_{\widetilde L^{\overline r}_t\dot B_{\infty,1}^{s+\frac{1}{\overline r}-1}}\big).
 $$
 Besides if $v=\nabla^\bot|\textnormal{D}|^{-1}\theta,$  then we have for all $s\geq1$ 
 $$
 \|\theta\|_{\widetilde L^r_{t}\dot{B}_{\infty,1}^{s+\frac{1}{r}}}\leq C e^{C\int_{0}^t(\| \nabla\theta(\tau)\|_{L^\infty}+\|\nabla v(\tau)\|_{L^\infty})d\tau}
 \big(  \|\theta^0\|_{\dot{B}_{\infty,1}^s}+\|f\|_{\widetilde L^{\overline r}_t\dot B_{\infty,1}^{s+\frac{1}{\overline r}-1}}\big).
 $$ 
 \end{Theo}
We use for the proof a new approach  based on
Lagrangian coordinates  combined with paradifferential calculus and a new commutator estimate. This idea has been recently used by the second author in \cite{Hmidi,Keraani}.

The rest of the paper is structured as follows. In section $2$ we review some basic results of Littlewood-Paley theory and we give some useful lemmas. Section $3$ deals with a new commutator estimate which is needed for the proof of Theorem \ref{Thm3}, done in section $4.$ Theorem \ref{Thm1} is proved in section $5.$

\section{Preliminaries}
 Throughout the paper, $C$ stands for a  constant which may be different in each occurrence. We shall sometimes use the notation $A\lesssim B$ instead of $A\leq CB$ and $A\approx B$ means that $A\lesssim B$ and $B\lesssim A.$
We denote by $\mathcal{ F}{f}$ the Fourier transform of $f$ and $[\eta]$ by the whole part of $\eta.$

One starts with recalling a traditional result that will be frequently used.
\begin{lem} $($Bernstein$)$ Let $f\in L^p(\RR^2),$ with  $1\leq p\leq\infty$ and $0<r<R.$
Then there exists a constant $C>0$ such that $\forall\;k\in\NN$ and $\lambda>0,$ we have:
$$
\begin{aligned}
&\mbox{supp}\,\mathcal{ F}f\subset B(0,\lambda r)
\Longrightarrow
\displaystyle\sup_{|\beta|=k}\|\partial^{\beta}f\|_{L^q}
\leq
C^k {\lambda}^{k+2(\frac{1}{p}-\frac{1}{q})}\|f\|_{L^p},
\\&
\mbox{supp}\,\mathcal{ F}f\subset\mathcal{C}(0,\lambda r,\lambda R)
\Longrightarrow
 C^{-k}{\lambda}^k\|f\|_{L^p}
 \leq
 \sup_{|\beta|=k}\|\partial ^{\beta}f\|_{L^p}
 \leq C^k{\lambda}^k\|f\|_{L^p}.
\end{aligned}
$$
These estimates hold true if we replace the derivation $\partial^\beta$ by $|\textnormal{D}|^{|\beta|}.$
\end{lem}

To define Besov spaces we need to recall the homogeneous Littlewood-Paley decomposition based on a dyadic unity partition. 
Let
$\varphi$ be a smooth function  supported in the ring $\mathcal{C}:=\{ \xi\in\RR^2,\frac{3}{4}\leq|\xi|\leq\frac{8}{3}\}$ and such that 
$$
\sum_{q\in\ZZ}\varphi(2^{-q}\xi)=1 \quad\hbox{for}\quad \xi\neq 0.
$$
Now,  for  $u\in{\mathcal S}'$ we set
  $$
\forall q\in\ZZ,\quad \Delta_qu=\varphi(2^{-q}\textnormal{D})u\hspace{1cm}\mbox{and}\hspace{1cm}
S_qu=\sum_{j\leq q-1}\Delta_{j}u.
 $$ 
We have the formal decomposition
$$
u=\sum_{q\in\ZZ}\Delta_q \,u,\quad\forall\,u\in {\mathcal {S}}'(\RR^2)/{\mathcal{P}}[\RR^2],
$$
where ${\mathcal{P}}[\RR^2]$ is the set of polynomials (see \cite{PE}).
Moreover, the Littlewood-Paley decomposition satisfies the
property of almost orthogonality:
\begin{equation}\label{Pres_orth}
\Delta_k\Delta_q u\equiv 0
\quad\mbox{if}\quad\vert k-q\vert\geq 2
\quad\mbox{and}\quad\Delta_k(S_{q-1}u\Delta_q u)
\equiv 0\quad\mbox{if}\quad\vert k-q\vert\geq 5.
\end{equation}
We recall now the definition of  Besov spaces. 
Let $(p,m)\in[1,+\infty]^2,$ $s\in\RR$ and $u\in{\mathcal S}',$ we set 
$$
\|u\|_{\dot B^s_{p,m}}:=\Big(2^{qs}\|\Delta_q u\|_{L^{p}}\Big)_{\ell ^{m}};\quad
\dot{\mathcal{B}}^s_{p,m}:=\Big\{u\in{\mathcal S}\;\big|\;
\|u\|_{\dot B^s_{p,m}}<\infty\Big\}\cdot 
$$
\begin{itemize}

\item 
For $s<\frac{2}{p}$ (or $s\leq\frac{2}{p}$ if $m=1$), we then define $\dot B^s_{p,m}$ as the completion of $\dot{\mathcal{B}}^s_{p,m}$ for
$\|\cdot\|_{\dot B^s_{p,m}}.$

\item 
If $k\in\NN$ and $\frac{2}{p}+k-1\leq s<\frac{2}{p}+k$ (or $s=\frac{2}{p}+k$ if $m=1$),
then $\dot B^s_{p,m}$ is defined as the subset of distributions $u\in{\mathcal S}'$ such that $\partial^\beta u\in\dot B^{s-k}_{p,m}$
whenever $|\beta|=k.$
\end{itemize}

Another characterization of the homogeneous Besov spaces that will be needed later is the following, see  for instance \cite{Triebel}. For $s\in]0,1[, p,m\in[1,\infty]$
\begin{equation}
\label{equivalence}
\Big(\int_{\RR^2}\frac{\|u(\cdot-x)-u(\cdot)\|^{m}_{L^p}}{|x|^{sm}}\frac{dx}{|x|^2}\Big)^{\frac{1}{m}}\approx\|u\|_{\dot{B}_{p,m}^s},
\end{equation}
with the usual modification if $m=\infty$.\\
In our next study we require two kinds of coupled space-time Besov spaces. The first one is defined by the following manner: for  $T>0$ \mbox{and $m\geq1,$} we denote by $L^r_{T}\dot B_{p,m}^s$ the set of all tempered distribution $u$ satisfying
$$
\|u\|_{L^r_{T}\dot B_{p,m}^s}:= \Big\|\Big( 2^{qs}
\|\Delta_q u\|_{L^{p}}\Big)_{\ell ^{m}}\Big\|_{L^r_{T}}<\infty.$$
The second mixed space is  $\widetilde L^r_{T}{\dot B_{p,m}^s}$ which is the set of  tempered distribution $u$ satisfying
 $$
 \|u\|_{ \widetilde L^r_{T}{\dot B_{p,m}^s}}:= \Big( 2^{qs}
\|\Delta_q u\|_{L^r_{T}L^{p}}\Big)_{\ell ^{m}}<\infty .$$
We can define by the same way the spaces $L^r_{T} B_{p,m}^s$ and $\widetilde L^r_{T}{ B_{p,m}^s}.$ 
The following embeddings  are a direct consequence of Minkowski's inequality. 

Let $s\in\RR,$ $r\geq1$ and $\big(p,m\big)\in[1,\infty]^2,$ then we have 
\begin{eqnarray}\label{lemm4}
L^r_{T}\dot B_{p,m}^s&\hookrightarrow&
\widetilde L^r_{T}{\dot B_{p,m}^s},\,\textnormal{if}\quad  m\geq r\quad\hbox{and}\\
\nonumber\widetilde L^r_{T}{\dot B_{p,m}^s}&
\hookrightarrow& L^r_{T}\dot B_{p,m}^s,\, \textnormal{if}\quad 
r\geq m.
\end{eqnarray}
The next lemma will be useful.
\begin{prop}\label{BERNSTEIN}
The following results hold true:

\begin{itemize}

\item 
$$
 \dot{B}_{p,m}^s\hookrightarrow\dot{B}_{p_{1},m_{1}}^{s-2(\frac{1}{p}-\frac{1}{p_{1}})},\quad\hbox{for}\quad p\leq p_{1}
 \quad\hbox{and}\quad m\leq m_{1}.
 $$
 
 \item
 Let $|\textnormal{D}|:=\sqrt{-\Delta}$ and $\sigma\in\RR,$ then the operator $|\textnormal{D}|^{\sigma}$ is an isomorphism from $\dot B^s_{p,m}$ to 
 $\dot B^{s-\sigma}_{p,m}.$
 \item
 Let $\gamma\in]0,1[,$ $s_1,$ $s_2\in\RR$ such that $s_1<s_2$ and $u\in\dot B^{s_1}_{p,\infty}\cap\dot B^{s_2}_{p,\infty},$ then
 $$
 \|u\|_{\dot B^{\gamma s_1+(1-\gamma)s_2}_{p,1}}
 \lesssim
 \|u\|_{\dot B^{s_1}_{p,\infty}}^{\gamma}
 \|u\|_{\dot B^{s_2}_{p,\infty}}^{1-\gamma}.
 $$
 
 \item
 For $s>0,$ $\dot B^s_{p,m}\cap L^\infty$ is an algebra.
 
 \end{itemize}
 \end{prop}

We now recall some  commutator estimates (see \cite{Ch1,D} and the references therein).
\begin{lem}\label{lemm12}
Let  $p,r\in[1,\infty], 1=\frac{1}{r}+\frac{1}{r'},  \rho_{1}<1, \rho_{2}<1$ and $v$ be a divergence free  vector field of $\RR^2.$ Assume in addition that
$$
\rho_{1}+\rho_{2}+2\min\{1,{2}/{p}\}>0\quad\hbox{and}
\quad\rho_{1}+{2}/{p}>0.
$$
Then we have
$$
\sum_{q\in\ZZ}2^{q(\frac{2}{p}+\rho_{1}+\rho_{2}-1)}
\big\|[\Delta_q, v\cdot \nabla ]u\big\|_{L^1_{t}L^{p}}
\lesssim  
\|v\|_{\widetilde L^r_{t}\dot B_{p,1}^{\frac{2}{p}+\rho_{1}}}\|u \|_{\widetilde L^{r'}_{t}\dot B_{p,1}^{\frac{2}{p}+\rho_{2}}}.
$$ 
Moreover we have for $s\in]-1,1[$
$$
\sum_{q\in\ZZ}2^{qs}\big\|[\Delta_q, v\cdot \nabla ]u\big\|_{L^{p}}\lesssim \|\nabla
v\|_{L^{\infty}}\|u\|_{\dot B_{p,1}^s}.
$$
If $v=\nabla^\bot|\textnormal{D}|^{-1}\theta,$ then the above estimate holds true for $s\geq 1$ if we \mbox{replace
$\|\nabla v\|_{L^\infty}$ }\mbox{by
$\|\nabla v\|_{L^\infty}+\|\nabla u\|_{L^\infty}.$}
\end{lem}
The following result is due to Vishik  \cite{v1}.

 \begin{lem}\label{l400} 
Let  $f$ be a function in  Schwartz class and $\psi$ a 
diffeomorphism   preserving  Lebesgue measure, then  for all
$p\in[1,+\infty]$ and for all  { $j,q\in\ZZ,$}
$$\|\Delta_j(\Delta_q f\circ\psi)\|_{L^p}\leq
C2^{-\vert j-q\vert}\|\nabla\psi ^{\epsilon(j,q)}\|_{L^{\infty}}\|\Delta_q
f\|_{L^p},$$
with
$$\epsilon(j,q)=\hbox{sign}(j-q).$$
\end{lem}
The following result is proved in \cite{C-C}.
\begin{prop}\label{maximum}
Let $v$ be a smooth divergence free vector field and $f$ be a smooth function. We assume that $\theta$ is a smooth solution of the equation
$$
\partial_{t}\theta+v\cdot\nabla \theta+\kappa |\textnormal{D}|^\alpha \theta=f, \quad\textnormal{with}\quad\kappa\geq 0\quad\textnormal{and}\quad
\alpha\in[0,2].
$$
Then for $p\in[1,+\infty]$ we have
$$
\|\theta(t)\|_{L^p}\leq\|\theta(0)\|_{L^p}+\int_{0}^t\|f(\tau)\|_{L^p}d\tau.
$$
\end{prop}

We can find a proof of the next proposition in \cite{Keraani}.
\begin{prop}\label{l:5}
Let $\mathcal{C}$ be a ring and $\alpha\in\RR_{+}.$ There exists a positive constant $C$ such
 that for any $p\in[1;+\infty],$ for any couple $(t,\lambda)$ of positive real numbers, we have
$$
\textnormal{supp}\, \mathcal{ F}u\subset\lambda\mathcal{C}\Rightarrow 
\|e^{-t|\textnormal{D}|^\alpha}u\|_{L^p}\leq Ce^{-C^{-1}t\lambda^{\alpha}}\|u\|_{L^p}.
$$
\end{prop}
\section{Commutator estimate}
The main result of this section is the following estimate that will play a  crucial role for the proof of Theorem \ref{Thm3}.
\begin{prop}\label{pr:1}
Let $v$ be a divergence free vector field belonging to  $L^1_{{\textnormal{loc}}}(\RR_{+};\textnormal{Lip}(\RR^2)).$ For  $q\in\ZZ$ 
we denote by $\psi_{q}$ the flow of the regularized vector field $S_{q-1}v.$ Then for $f\in\dot{B}_{\infty,\infty}^1$ and for $q\in\ZZ$ we have
\begin{eqnarray*}
\big{\|}|\textnormal{D}|(\Delta_{q}f\circ\psi_{q})-(|\textnormal{D}| \Delta_{q}f)\circ\psi_{q}\big{\|}_{L^\infty}&\leq&
 Ce^{CV(t)}V^{\frac{1}{2}}(t)2^q\|\Delta_{q}f\|_{L^\infty},
\end{eqnarray*}
where $V(t)=\|\nabla v\|_{L^1_{t}L^\infty(\RR^2)}$ and $C$   an absolute constant.
\end{prop}
{\bf{Proof}.} 
We set $f_{q}:=\Delta_{q}f,$ then it is obvious that 
\begin{eqnarray*}
|\textnormal{D}|(f_{q}\circ\psi_{q})-(|\textnormal{D}| f_{q})\circ\psi_{q}&=&
|\textnormal{D}|^{\frac{1}{2}}\{(|\hbox{D}|^{\frac{1}{2}}f_{q})\circ\psi_{q}\}
-\{|\textnormal{D}|^{\frac{1}{2}}(|\hbox{D}|^{\frac{1}{2}} f_{q})\}\circ\psi_{q}\\
&+&|\textnormal{D}|^{\frac{1}{2}}\big\{|\hbox{D}|^{\frac{1}{2}}(f_{q}\circ\psi_{q})-(|\textnormal{D}|^{\frac{1}{2}} f_{q})\circ\psi_{q}\big\}\\
&:=& \hbox{I}+\hbox{II}.
\end{eqnarray*}
For the first term we apply Proposition 3.1 \cite{Keraani}, with $\alpha=\frac{1}{2}$ and 
$F_{q}=|\hbox{D}|^{\frac{1}{2}}f_{q},$ yielding
\begin{eqnarray*}
\|\hbox{I}\|_{L^\infty}&\leq& Ce^{CV(t)}(e^{CV(t)}-1)\|F_{q}\|_{\dot{B}_{\infty,1}^{\frac{1}{2}}}\\
&\leq&
Ce^{CV(t)}(e^{CV(t)}-1)2^q\|f_{q}\|_{L^\infty}\\
&\lesssim&e^{CV(t)}V^{\frac{1}{2}}(t)2^q\|f_{q}\|_{L^\infty}.
\end{eqnarray*}

 For the second term we use the following formula for the fractional Laplacian 
  $$
|\textnormal{D}|^{\frac{1}{2}} f(x)=C\,\int_{\RR^2}\frac{f(x)-f(y)}{\vert x-y \vert^{\frac{5}{2}}}dy.
$$ 
Since the flow $\psi_{q}$ preserves Lebesgue measure then we get easily
 \begin{eqnarray*}
 |\hbox{D}|^{\frac{1}{2}}(f_{q}\circ\psi_{q})(x)-(|\hbox{D}|^{\frac{1}{2}} f_{q})\circ\psi_{q}(x)
&=& C\int_{\RR^2}\frac{f_{q}(\psi_{q}(x))-f_{q}(\psi_{q}(y))}{\vert x-y \vert^{\frac{5}{2}}}\\
&&\times\Big( 1-\frac{\vert x-y \vert^{\frac{5}{2}}}{\vert \psi_{q}(x)-\psi_{q}(y) \vert^{\frac{5}{2}}} \Big) dy.
\end{eqnarray*}
We denote $g_{q}(x)=f_{q}(\psi_{q}(x))$ and we put $h=x-y:$
\begin{eqnarray*}
 |\hbox{D}|^{\frac{1}{2}}(f_{q}\circ\psi_{q})(x)-(|\hbox{D}|^{\frac{1}{2}} f_{q})\circ\psi_{q}(x)
&=& C\int_{\RR^2}\frac{g_{q}(x)-g_{q}(x-h)}{\vert h \vert^{\frac{5}{2}}}\bar{\psi}_{q}(x,h)dh,\end{eqnarray*}
with
$$
\bar{\psi}_{q}(x,h)= 1-\frac{\vert h \vert^{\frac{5}{2}}}{\vert \psi_{q}(x)-\psi_{q}(x-h) \vert^{\frac{5}{2}}}\cdot
$$
It follows from law products and the embedding $\dot{B}_{\infty,1}^0\hookrightarrow L^\infty$
\begin{eqnarray*}
\|\hbox{II}\|_{L^\infty}
&\leq& 
\big{\|}|\hbox{D}|^{\frac{1}{2}}(f_{q}\circ\psi_{q})-(|\hbox{D}|^{\frac{1}{2}} f_{q})\circ\psi_{q}\big{\|}_{\dot{B}_{\infty,1}^{\frac{1}{2}}}\\
&\leq&
 C\|\bar{\psi}_{q}\|_{L^\infty(\RR^{4})}\int_{\RR^2}|h|^{-\frac{5}{2}}
 {\|g_{q}(\cdot)-g_{q}(\cdot-h)\|_{\dot{B}_{\infty,1}^{\frac{1}{2}}}}{dh}\\
&+&
C\sup_{h\in\RR^2}\|\bar{\psi}_{q}(\cdot,h)\|_{\dot{B}_{\infty,1}^{\frac{1}{2}}}\int_{\RR^2}|h|^{-\frac{5}{2}}\|g_{q}(\cdot)-g_{q}(\cdot-h)\|_{L^\infty}dh\\
&=& 
J_{q}^1+J_{q}^2. 
\end{eqnarray*}
We intend to estimate  $J_{q}^1.$
It is plain from
 Mean Value Theorem that
 $$
\frac{1}{\|\nabla \psi\|_{L^\infty}^{\frac{5}{2}}}\leq\frac{|h|^{\frac{5}{2}}}{|\psi(x)-\psi(x-h)|^{\frac{5}{2}}}\leq \|\nabla\psi^{-1}\|_{L^\infty}^{\frac{5}{2}},
$$ 
 
which gives  easily the inequality
$$
{\|} \bar{\psi}_{q} {\|}_{L^\infty(\RR^{4})}
 \leq \max\Big(\vert 1-\|\nabla \psi_{q}^{-1} \|_{{L^\infty}}^{\frac{5}{2}}\vert;
 \vert 1-{\|\nabla \psi_{q} \|_{{L^\infty}}^{-\frac{5}{2}}}\vert \Big). 
  $$
On the other hand we have the classical estimates
\begin{equation}\label{Q0}
\begin{aligned}
e^{-C\|S_{q-1}\nabla v\|_{L^1_{t}L^\infty}}
&\leq
\|\nabla\psi_{q}^{\mp1}\|_{L^\infty}
\leq 
e^{C\|S_{q-1}\nabla v\|_{L^1_{t}L^\infty}}
\\&
\hbox{and}\quad
\|S_{q-1}\nabla v\|_{L^1_{t}L^\infty}\leq CV(t).
\end{aligned}
\end{equation}
We thus get
\begin{equation}\label{Q1} 
{\|} \bar{\psi}_{q} {\|}_{L^\infty(\RR^{4})}\leq Ce^{CV(t)}(e^{CV(t)}-1).\end{equation}
Using the definition of Besov spaces and the commutation of $\Delta_{j}$ with translation operators one finds
$$
\begin{aligned}
\int_{\RR^2}|h|^{-\frac{5}{2}}\|g_{q}(\cdot)&-g_{q}(\cdot-h)\|_{\dot{B}_{\infty,1}^{\frac{1}{2}}}{dh}
\\&
\leq \sum_{j}2^{\frac{1}{2}j}\int_{\RR^2}|h|^{-\frac{1}{2}}\|\Delta_{j}g_{q}(\cdot)-(\Delta_{j}g_{q})(\cdot-h)\|_{L^\infty}{\frac{dh}{|h|^2}}\cdot
\end{aligned}
$$
Applying the characterization of Besov spaces (\ref{equivalence}) yields
\begin{eqnarray*}
\int_{\RR^2}|h|^{-\frac{5}{2}}\|g_{q}(\cdot)-g_{q}(\cdot-h)\|_{\dot{B}_{\infty,1}^{\frac{1}{2}}}{dh}
&\leq&
C\sum_{j}2^{\frac{1}{2}j}\|\Delta_{j} g_{q}\|_{\dot B^{\frac{1}{2}}_{\infty,1}}\\
&\leq& C\sum_{|j-k|\leq1}2^{j\frac{1}{2}}2^{\frac{1}{2}k}\|\Delta_{j}\Delta_k g_{q}\|_{L^\infty}\\
&\leq& C\|g_{q}\|_{\dot{B}_{\infty,1}^1}.
\end{eqnarray*}
Now  we use the following interpolation estimate
\begin{eqnarray*}
\|g_{q}\|_{\dot{B}_{\infty,1}^1}&\lesssim&\|g_{q}\|_{L^\infty}^{\frac{1}{2}}\|\Delta g_{q}\|_{L^\infty}^{\frac{1}{2}}\\
&\lesssim& \|f_{q}\|_{L^\infty}^{\frac{1}{2}}\|\Delta g_{q}\|_{L^\infty}^{\frac{1}{2}}.
\end{eqnarray*}
It is easy to check from Leibnitz rule that
$$
\Delta g_{q}=\Delta (f_{q}\circ\psi_{q})=\sum_{i=1}^d\langle (\nabla^2 f_{q})\circ\psi_{q}\cdot\partial_{i}\psi_{q},\partial_{i}\psi_{q}\rangle+(\nabla f_{q})\circ\psi_{q}\cdot\Delta\psi_{q}.
$$
Applying Bernstein inequality we get
$$
\|\Delta g_{q}\|_{L^\infty}\lesssim e^{CV(t)} 2^{2q}\|f_{q}\|_{L^\infty}+2^q\|f_{q}\|_{L^\infty}\|\Delta\psi_{q}\|_{L^\infty}.
$$
The derivative of the flow equation with respecyt to $x$ and the use of Gronwall and Bernstein inequalities  gives
\begin{eqnarray}\label{Q5}
\nonumber\|\nabla^2\psi_{q}(t)\|_{L^\infty}&\lesssim& e^{CV(t)}\int_{0}^t\|\nabla^2S_{q-1}v(\tau)\|_{L^\infty}d\tau\\
&\lesssim& e^{CV(t)}2^q.
\end{eqnarray}
Combining both last estimates we obtain
\begin{equation}\label{Q2}
\|\Delta g_{q}\|_{L^\infty}\lesssim e^{CV(t)} 2^{2q}\|f_{q}\|_{L^\infty}.
\end{equation}
Putting together (\ref{Q1}) and (\ref{Q2}) we conclude that
$$
\|J_{q}^1(t)\|_{L^\infty}\lesssim e^{CV(t)}(e^{CV(t)}-1) 2^{q}\|f_{q}\|_{L^\infty}.
$$
Let us now turn to the second term $J_{q}^2.$ The integral term can be estimated from (\ref{equivalence}) as follows
\begin{eqnarray*}
\int_{\RR^2}|h|^{-\frac{5}{2}}\|g_{q}(\cdot)-g_{q}(\cdot-h)\|_{L^\infty}dh\lesssim \|g_{q}\|_{\dot{B}_{\infty,1}^{\frac{1}{2}}}.
\end{eqnarray*}
According to classical composition result we write   
 \begin{eqnarray}\label{Q11}
 \nonumber\|g_{q}(t)\|_{\dot{B}_{\infty,1}^{\frac{1}{2}}}&\lesssim&\|\nabla\psi_{q} \|_{L^\infty}^{\frac{1}{2}}\|f_{q}\|_{\dot{B}_{\infty,1}^{\frac{1}{2}}}\\
 &\lesssim& e^{CV(t)}2^{q\frac{1}{2}}\|f_{q}\|_{L^\infty}.
 \end{eqnarray}
 In order to estimate $\bar{\psi}_{q}$ we use the interpolation inequality
 \begin{equation*}
\|\bar{\psi}_{q}(\cdot,h)\|_{\dot{B}_{\infty,1}^{\frac{1}{2}}}\lesssim
\|\bar{\psi}_{q}(\cdot,h)\|_{L^\infty}^{\frac{1}{2}}\|\nabla_{x}\bar{\psi}_{q}(\cdot,h)\|_{L^\infty}^{\frac{1}{2}}.
\end{equation*}
This leads in view  of (\ref{Q1}) to
 \begin{equation}\label{Q6}\|\bar{\psi}_{q}(\cdot,h)\|_{\dot{B}_{\infty,1}^{\frac{1}{2}}}\leq C
e^{CV(t)}(e^{CV(t)}-1)^{\frac{1}{2}}\|\nabla_{x}\bar{\psi}_{q}(\cdot,h)\|_{L^\infty}^{\frac{1}{2}}.
\end{equation}
The derivative of $\bar{\psi}_{q}$ with respect to $x$ yields
\begin{eqnarray*}
|\nabla_{x}\bar{\psi}_{q}(x,h)|&\lesssim &\frac{|h|^{\frac{7}{2}}}
{|\psi_{q}(x)-\psi_{q}(x-h)|^{\frac{7}{2}}}\frac{|\nabla_{x}\psi_{q}(x)-\nabla_{x}\psi_{q}(x-h)|}{|h|}\\
&\lesssim&\|\nabla\psi_{q}^{-1}\|_{L^\infty}^{\frac{7}{2}}\|\nabla^2\psi_{q}\|_{L^\infty}.
\end{eqnarray*}
Combining (\ref{Q0}) and (\ref{Q5})
 we obtain
 \begin{equation}\label{Q8}
 \|\nabla_{x}\bar{\psi}_{q}(t)\|_{L^\infty(\RR^{4})}\lesssim e^{CV(t)}2^q.
 \end{equation}
 Plugging (\ref{Q8}) into (\ref{Q6}) we find
 \begin{equation}\label{Q10}
 \|\bar{\psi}_{q}(\cdot,h)\|_{\dot{B}_{\infty,1}^{\frac{1}{2}}}\lesssim
e^{CV(t)}V^{\frac{1}{2}}(t)2^{\frac{q}{2}}.
\end{equation}
 Thus  we deduce from (\ref{Q10}) and (\ref{Q11}) that
 $$
 \|J_{q}^2(t)\|_{L^\infty}\leq Ce^{CV(t)}V^{\frac{1}{2}}(t)2^{q}\|f_{q}(t)\|_{L^\infty}. $$
 This achieves the proof.
$\hfill$$\square$

\section{Proof of Theorem \ref{Thm3}}
The Fourier localized function   $\theta_{q}:=\Delta_{q}\theta$ satisfies
 
 \begin{equation}\label{eq:1}
 \partial_{t}\theta_{q}+S_{q-1}v\cdot\nabla\theta_{q}+|\hbox{D}|\theta_{q}=-[\Delta_{q},v\cdot\nabla]\theta+(S_{q-1}v-v)\cdot\nabla\theta_{q}
 +f_{q}:=\mathcal{R}_{q}.
 \end{equation}

 Let  $\psi_{q}$ denote the flow of the  velocity $S_{q-1}v$ and  set 
 $$
 \bar\theta_{q}(t,x)=\theta_{q}(t,\psi_{q}(t,x))\quad\hbox{and}\quad \bar {\mathcal R}_{q}(t,x)=\mathcal{R}_{q}(t,\psi_{q}(t,x)).
 $$
Since $\psi_q$ is an homeomophism, then
 \begin{equation}
 \label{eq:5}
\|\bar{\mathcal{R}}_{q}\|_{L^\infty}\leq \|[\Delta_{q},v\cdot\nabla]\theta\|_{L^\infty}+\|(S_{q-1}v-v)\cdot\nabla\theta_{q}\|_{L^\infty}
+\|f_{q}\|_{L^\infty}.
\end{equation}
It is not hard to check that the function $\bar\theta_{q}$ satisfies 
 \begin{equation}
 \label{T1}
 \partial_{t}\bar\theta_{q}+|\hbox{D}|\bar\theta_{q}=
 |\hbox{D}|(\theta_{q}\circ\psi_{q})-(|\hbox{D}|\theta_{q})\circ\psi_{q}+\bar {\mathcal R}_{q}:=\bar{\mathcal R}_{q}^1.
 \end{equation}
From Proposition \ref{pr:1} we find that for  $q\in\ZZ$
\begin{equation}
\label{T:5}
\| |\hbox{D}|(\theta_{q}\circ\psi_{q})-(|\hbox{D}|\theta_{q})\circ\psi_{q}\|_{L^\infty}
 \lesssim e^{CV(t)} V^{\frac{1}{2}}(t)2^q\|\theta_{q}(t)\|_{L^\infty}. \end{equation}
where $V(t):=\|\nabla v\|_{L^1_{t}L^\infty}.$ Putting together (\ref{eq:5}) and (\ref{T:5}) yields
\begin{eqnarray*}
\|\bar{\mathcal R}_{q}^1(t)\|_{L^\infty}\lesssim \| f_{q}(t)\|_{L^\infty}&+&\|(S_{q-1}v-v)\cdot\nabla\theta_{q}\|_{L^\infty}
+\|[\Delta_{q},v\cdot\nabla]\theta\|_{L^\infty}\\
&+&e^{CV(t)} V^{\frac{1}{2}}(t)2^{q}\|\theta_{q}(t)\|_{L^\infty}.
\end{eqnarray*}
Applying the operator $\Delta_{j}$  to the equation (\ref{T1}) and using Proposition \ref{l:5}
\begin{eqnarray} 
 \|\Delta_{j}\bar\theta_{q}(t)\|_{L^\infty}
 &\lesssim& 
 e^{-ct2^{j}}\|\Delta_{j}\theta_{q}^0\|_{L^\infty}+\int_{0}^te^{-c (t-\tau)2^{j}}\| f_{q}(\tau)\|_{L^\infty}d\tau \\
\nonumber&+& e^{CV(t)}V^{\frac{1}{2}}(t)2^{q}
 \int_{0}^te^{-c (t-\tau)2^{j}}\|\theta_{q}(\tau)\|_{L^\infty}d\tau
 \\\nonumber &+&
 \int_{0}^te^{-c (t-\tau)2^{j}}\|[\Delta_{q},v\cdot\nabla]\theta(\tau)\|_{L^\infty}d\tau
\\&+&
\int_{0}^te^{-c (t-\tau)2^{j}}\|(S_{q-1}v-v)\cdot\nabla\theta_{q}(\tau)\|_{L^\infty}d\tau.  
\end{eqnarray}
Integrating this estimate with respect to the time  and using Young inequality 
\begin{eqnarray}\label{eres}
\nonumber \|\Delta_{j}\bar\theta_{q}\|_{L^r_{t}L^\infty}
&\lesssim& 
2^{-\frac{j}{r}}(1-e^{-crt2^{j}})^{\frac{1}{r}}\|\Delta_{j}\theta_{q}^0\|_{L^\infty}
+2^{-j(1+\frac{1}{r}-\frac{1}{\overline r})}\|f_{q}\|_{L^{\overline r}_{t}L^\infty}
\\\nonumber&+&
e^{CV(t)}V^{\frac{1}{2}}(t)2^{(q-j)}\|\theta_{q}\|_{L^r_{t}L^\infty}\\
\nonumber&+&2^{-\frac{j}{r}}\int_{0}^t \|[\Delta_{q},v\cdot\nabla]\theta(\tau)\|_{L^\infty}d\tau
\\&+&
2^{-\frac{j}{r}}\int_{0}^t\|(S_{q-1}v-v)\cdot\nabla\theta_{q}(\tau)\|_{L^\infty}d\tau.
\end{eqnarray}
Since the flow $\psi$ is an homeomorphism then one writes
\begin{eqnarray*}&&2^{q(s+\frac{1}{r})}\|\theta_{q}\|_{L^r_{t}L^\infty}=2^{q(s+\frac{1}{r})}
\|\bar\theta_{q}\|_{L^r_{t}L^\infty}
\\&\leq &
2^{q(s+\frac{1}{r})}
\Big(\sum_{\vert j-q \vert >N}\|\Delta_{j}\bar\theta_{q}\|_{L^r_{t}L^\infty}+
\sum_{\vert j-q \vert\leq N}\|\Delta_{j}\bar\theta_{q}\|_{L^r_{t}L^\infty}\Big)
\\
&:=&
\textnormal{I}+\textnormal{II}.
\end{eqnarray*}
To estimate the term $\textnormal I$ we make appeal to  Lemma \ref{l400} 
\begin{eqnarray*}
\|\Delta_{j}\bar\theta_{q}\|_{L^r_{t}L^\infty}
&\lesssim& 
2^{-\vert q-j\vert}
e^{C\int_{0}^t
\|\nabla v(\tau)\|_{L^\infty}d\tau}\|\theta_{q}\|_{L^r_{t}L^\infty}\\
&\leq&
C2^{-\vert q-j\vert}
e^{CV(t)}\|\theta_{q}\|_{L^r_{t}L^\infty}.
\end{eqnarray*}
Therefore we get
\begin{equation}\label{e:11}
\textnormal{I}
\leq
C2^{-N}e^{CV(t)}2^{q(s+\frac{1}{r})}\|\theta_{q}\|_{L^r_{t}L^\infty}.
\end{equation}
In order to bound the second term $\textnormal{II}$ we use (\ref{eres}) 
\begin{eqnarray}\label{Meth1}
\nonumber \textnormal{II}
&\lesssim& 
(1-e^{-crt2^{q}})^{\frac{1}{r}}2^{qs}\|\theta_{q}^0\|_{L^\infty}+2^{N(\frac{1}{r}+1-\frac{1}{\overline r})}2^{q(s+\frac{1}{\overline r}-1)}
\|f_{q}\|_{L^{\overline r}_{t}L^\infty}
\\\nonumber&+&
2^{N}e^{CV(t)}V^{\frac{1}{2}}(t)2^{q(s+\frac{1}{r})}\|\theta_{q}\|_{L^r_{t}L^\infty}\\
\nonumber&+&
2^{\frac{N}{r}}2^{qs}\int_{0}^t \|[\Delta_{q},v\cdot\nabla]\theta(\tau)\|_{L^\infty}d\tau\\
&+&
2^{\frac{N}{r}}2^{qs}\int_{0}^t\|(S_{q-1}v-v)\cdot\nabla\theta_{q}(\tau)\|_{L^\infty}d\tau.
\end{eqnarray}
Denote  $Z_{q}^r(t):=2^{q(s+\frac{1}{r})}\|\theta_{q}\|_{L^r_{t}L^\infty},$ then we obtain  in view of  (\ref{e:11}) and (\ref{Meth1}) 
\begin{eqnarray*}
Z_{q}^r(t)
&\leq &
C(1-e^{-crt2^{q}})^{\frac{1}{r}}2^{qs}\|\theta_{q}^0\|_{L^\infty}+C2^{N(\frac{1}{r}+1-\frac{1}{\overline r})}2^{q(s+\frac{1}{\overline r}-1)}
\|f_{q}\|_{L^{\overline r}_{t}L^\infty}
\\&+&
C\big\{2^{N}e^{CV(t)}
V^{\frac{1}{2}}(t)+2^{-N}e^{CV(t)}\big\}Z_{q}^r(t)
\\&+&
C2^{\frac{N}{r}}2^{qs}\int_{0}^t \|[\Delta_{q},v\cdot\nabla]\theta(\tau)\|_{L^\infty}d\tau\\
&+&
C2^{\frac{N}{r}}2^{qs}\int_{0}^t\|(S_{q-1}v-v)\cdot\nabla\theta_{q}(\tau)\|_{L^\infty}d\tau.
\end{eqnarray*}
It is easy to check  the existence of two absolute constants $N$ and $C_{0}$ such that
$$
V(t)\leq C_{0}\Rightarrow C2^{-N}e^{CV(t)}+C2^{N}e^{CV(t)}V^{\frac{1}{2}}(t)\leq \frac{1}{2}\cdot
$$ 
Thus we obtain under this condition
\begin{eqnarray}
\label{mahma00}
\nonumber Z_{q}^r(t)
&\lesssim&
 (1-e^{-crt2^{q}})^{\frac{1}{r}}2^{qs}\|\theta_{q}^0\|_{L^\infty}+2^{q(s+\frac{1}{\overline r}-1)}
\|f_{q}\|_{L^{\overline r}_{t}L^\infty}
 \\&+&
2^{qs}\int_{0}^t\Big( \|[\Delta_{q},v\cdot\nabla]\theta(\tau)\|_{L^\infty}+\|(S_{q-1}v-v)\cdot\nabla\theta_{q}\|_{L^\infty}\Big)d\tau.
\end{eqnarray} 
Summing over $q$  and using Lemma \ref{lemm12} leads  for $V(t)\leq C_{0},$
\begin{eqnarray}\label{since}
\nonumber \|\theta\|_{\widetilde L^r_{t}\dot{B}_{\infty,1}^{s+\frac{1}{r}}}
 &\lesssim&
 \|\theta^0\|_{\dot{B}_{\infty,1}^s}+\|f\|_{\widetilde L^{\overline r}_{t}\dot B_{\infty,1}^{s+\frac{1}{\overline r}-1}}+
 \int_{0}^t\| \nabla v(\tau)\|_{L^\infty}
\|\theta(\tau)\|_{\dot{B}_{\infty,1}^s}d\tau\\
&\lesssim&\|\theta^0\|_{\dot{B}_{\infty,1}^s}+\|f\|_{\widetilde L^{\overline r}_{t}\dot B_{\infty,1}^{s+\frac{1}{\overline r}-1}}+C_{0}
 \|\theta\|_{L^\infty_{t}\dot{B}_{\infty,1}^s}.
 \end{eqnarray}
Let us show how to conclude the proof in the case  of $r=\infty.$ If $C_{0}$ is sufficiently small then we obtain from (\ref{since}) the desired estimate:
\begin{eqnarray}\label{1er_etape}
 \|\theta\|_{\widetilde L^\infty_{t}\dot{B}_{\infty,1}^{s}}
 \lesssim
\|\theta^0\|_{\dot{B}_{\infty,1}^s}+\|f\|_{\widetilde L^{\overline r}_{t}\dot B_{\infty,1}^{s+\frac{1}{\overline r}-1}}.
\end{eqnarray}
Now for an arbitrary positive time  $T$ we take a partition $(T_{i})_{i=0}^M$  of $[0,T],$ 
such \mbox{that $\displaystyle \int_{T_{i}}^{T_{i+1}}\|\nabla v(\tau)\|_{L^\infty}d\tau\approx C_{0}.$}
We can proceed analogously to the above calculus and obtain 
$$
 \|\theta\|_{\widetilde L^\infty_{[T_{i},T_{i+1}]}\dot{B}_{\infty,1}^{s}}
\lesssim
 \|\theta(T_{i})\|_{\dot{B}_{\infty,1}^s}+\|f\|_{\widetilde L^{\overline r}_{[T_i,T_{i+1}]}\dot B_{\infty,1}^{s+\frac{1}{\overline r}-1}}.
$$
An iteration argument leads to
 $$
 \|\theta\|_{\widetilde L^\infty_{[T_{i},T_{i+1}]}\dot{B}_{\infty,1}^{s}}
\leq
C^{i+1}
\Big( \|\theta^0\|_{\dot{B}_{\infty,1}^s}+\|f\|_{\widetilde L^{\overline r}_{[0,T_{i+1}]}\dot B_{\infty,1}^{s+\frac{1}{\overline r}-1}}\Big).
$$
The triangle inequality and the fact that $C_{0}M\simeq 1+V(t)$ give
\begin{equation}\label{neige1}
\|\theta\|_{\widetilde L^\infty_{T}\dot{B}_{\infty,1}^{s}}
\leq 
Ce^{C \int_{0}^T\| \nabla v(\tau)\|_{L^\infty}}
\Big( \|\theta^0\|_{\dot{B}_{\infty,1}^s}+\|f\|_{\widetilde L^{\overline r}_{T}\dot B_{\infty,1}^{s+\frac{1}{\overline r}-1}}\Big).
\end{equation}
Let us now turn to the case of finite  $r.$ Combining (\ref{since}) and (\ref{1er_etape}) we obtain under the assumption $V(t)\leq C_{0}$
\begin{equation}\label{EDF}
 \|\theta\|_{\widetilde L^r_{t}\dot{B}_{\infty,1}^{s+\frac{1}{r}}}
 \lesssim
\|\theta^0\|_{\dot{B}_{\infty,1}^s}+\|f\|_{\widetilde L^{\overline r}_{t}\dot B_{\infty,1}^{s+\frac{1}{\overline r}-1}}.
 \end{equation}
This gives the result for a short time and  as for the case $r=\infty$ we obtain the required global estimate. \\
Concerning the last  estimate of Theorem \ref{Thm3}, we use in the commutator term \mbox{of (\ref{mahma00})} the last part of Lemma \ref{lemm12}.
$\hfill$$\square$

 \section{Proof of Theorem {\ref{Thm1}}}
  The proof is divided into two parts: in the first one we construct  local unique solution and we give a criteria of global existence. However we discuss in the second part the global existence by reproducing the same idea of \cite{volberg}. 
  \subsection{Local existence.}
We aim to prove the following result.

\begin{prop}\label{N67}
Given any $\theta^0\in\dot B^0_{\infty,1},$ there is  $T>0$ such that the $\QG$ equation has a unique solution $\theta$ with
$$
\theta\in\widetilde L^\infty_T\dot B^0_{\infty,1}\cap L^1_T\dot B^1_{\infty,1}.
$$
Moreover for all $\beta\in\RR_+,$ we have $t^\beta\theta\in\widetilde L^\infty_T\dot B^\beta_{\infty,1}.$
 \end{prop}
\begin{proof}
The existence is based on Theorem \ref{Thm3} and an iterative method.
We denote  $\theta_0(t,x):=e^{-t|\textnormal{D}|
}\theta^0(x),$ $v_0:=(-R_2\theta_0,R_1\theta_0)$ and $\theta_{n+1}$ the solution of the linear system
$$
\left\lbrace
\begin{array}{l}\partial_{t}\theta_{n+1}+v_n\cdot\nabla\theta_{n+1}
+|\hbox{D}|\theta_{n+1}=0,\\
v_{n}=(-R_{2} \theta_{n},R_{1}\theta_{n}),\\
{\theta_{n+1}}_{|t=0}=\theta^0.
\end{array}
\right.
$$
Since $\theta_0\in L^1(\RR_+;\,\dot B^1_{\infty,1})$ and from the continuity of Riesz transforms in the homogeneous Besov spaces we find
$v_0\in L^1(\RR_+;\,\dot B^1_{\infty,1}).$ Thus by iteration and thanks to  Theorem \ref{Thm3}, one deduces  that $\forall n\in\NN,$ $$
\theta_n\in\widetilde L^\infty(\RR_+;\,\dot B^0_{\infty,1})\cap L^1(\RR_+;\,\dot B^1_{\infty,1}).
$$
{\it Step 1: uniform bounds.}

Now we intend to obtain  uniform bounds, with respect to the parameter $n,$ for some $T>0$ independent of $n.$

By (\ref{mahma00}), we have for all $T\geq0$ such that \begin{equation}\label{recurrence}
\int_0^T\|\theta_n(\tau)\|_{\dot B^1_{\infty,1}}d\tau
\leq
C_1(:=CC_{0})
\end{equation}
the following estimate
$$
\begin{aligned}
\|\theta_{n+1}\|_{\widetilde L^2_T\dot B^{\frac{1}{2}}_{\infty,1}}
&+
\|\theta_{n+1}\|_{L^1_T\dot B^1_{\infty,1}}
\lesssim
\sum_{q\in\ZZ}(1-e^{-cT2^q})^{\frac{1}{2}}\|\Delta_q\theta^0\|_{L^\infty}
\\&
+
\sum_{q\in\ZZ}\int_0^T\|[\Delta_q,v_n\cdot\nabla]\theta_{n+1}(\tau)
\|_{L^\infty}d\tau
\\
&
+
\sum_{q\in\ZZ}\int_0^T\|(S_{q-1}v_n-v_n)\cdot\nabla\Delta_q\theta_{n+1}(\tau)\|_{L^\infty}d\tau.
\end{aligned}
$$
Since ${\mathop{\rm div}}\, v_n=0,$ then Lemma \ref{lemm12}  combined with the continuity of Riesz transforms 
gives
$$
\begin{aligned}
\sum_{q\in\ZZ}\int_0^T\|[\Delta_q,v_n\cdot\nabla]\theta_{n+1}(\tau)\|_{L^\infty}d\tau
&\lesssim
\|v_n\|_{\widetilde L^2_T\dot B^{\frac{1}{2}}_{\infty,\infty}}
\|\theta_{n+1}\|_{\widetilde L^2_T\dot B^{\frac{1}{2}}_{\infty,1}}
\\&
\lesssim
\|\theta_n\|_{\widetilde L^2_T\dot B^{\frac{1}{2}}_{\infty,\infty}}
\|\theta_{n+1}\|_{\widetilde L^2_T\dot B^{\frac{1}{2}}_{\infty,1}}.
\end{aligned}
$$
We deduce from H\"older and Young inequalities
$$
\begin{aligned}
\sum_{q\in\ZZ}\int_0^t\|(S_qv_n-&v_n)\cdot\nabla\Delta_q\theta_{n+1}\|_{L^\infty}d\tau
\lesssim
\sum_{q\in\ZZ}2^q\|\Delta_q\theta_{n+1}\|_{L^2_tL^\infty}\|S_qv_n-v_n\|_{L^2_tL^\infty}
\\&
\lesssim
\sum_{q\in\ZZ}2^{\frac{1}{2}q}\|\Delta_q\theta_{n+1}\|_{L^2_tL^\infty}
\sum_{k\geq q}2^{\frac{1}{2}(q-k)}2^{\frac{1}{2}k}
\|\Delta_kv_n\|_{L^2_tL^\infty}
\\&
\lesssim
\|\theta_{n+1}\|_{\widetilde L^2_t\dot B^{\frac{1}{2}}_{\infty,1}}
\|\theta_n\|_{\widetilde L^2_t\dot B^{\frac{1}{2}}_{\infty,\infty}}.
\end{aligned}
$$
Therefore we obtain from the above inequalities 
$$
\begin{aligned}
\|\theta_{n+1}\|_{\widetilde L^2_t\dot B^{\frac{1}{2}}_{\infty,1}}
+
\|\theta_{n+1}\|_{L^1_t\dot B^1_{\infty,1}}
&\lesssim
\sum_{q\in\ZZ}(1-e^{-ct2^q})^{\frac{1}{2}}\|\Delta_q\theta^0\|_{L^\infty}
\\
&
+
\|\theta_{n+1}\|_{\widetilde L^2_t\dot B^{\frac{1}{2}}_{\infty,1}}
\|\theta_n\|_{\widetilde L^2_t\dot B^{\frac{1}{2}}_{\infty,1}}.
\end{aligned}
$$
Thus there exists an absolute constant  $\varepsilon_0>0$ such that, if
\begin{equation}\label{temps_existence}
\sum_{q\in\ZZ}(1-e^{-cT2^q})^{\frac{1}{2}}\|\Delta_q\theta^0\|_{L^\infty}
\leq
\varepsilon_0,
\end{equation}
then
\begin{equation}\label{condition_petitesse}
\|\theta_{n+1}\|_{\widetilde L^2_T\dot B^{\frac{1}{2}}_{\infty,1}}
+
\|\theta_{n+1}\|_{L^1_T\dot B^1_{\infty,1}}
\leq
2\varepsilon_0.
\end{equation}
The existence of $T>0$ is due to Lebesgue theorem.

Hence, by using the estimate
$$
\int_{0}^T\|\nabla v(\tau)\|_{L^\infty}d\tau\lesssim \int_{0}^T\|\theta(\tau)\|_{\dot{B}_{\infty,1}^1}d\tau$$
 and Theorem \ref{Thm3} we obtain
$$
\|\theta_{n+1}\|_{\widetilde L^\infty_T\dot B^0_{\infty,1}}
\lesssim\|\theta^0\|_{\dot B^0_{\infty,1}}.
$$
Thus we prove that the sequence $(v_n,\theta_n)_{n\in\NN}$ is uniformly bounded in the space 
$\widetilde L^\infty_T\dot B^0_{\infty,1}\cap L^1_T\dot B^1_{\infty,1}.$
\\

{\it Step 2: strong convergence}.\\
We will prove that the sequence $(v_n,\theta_n)$ is of Cauchy in
$\widetilde L^\infty_T\dot B^0_{\infty,1}.$

Let $(n,m)\in\NN^2,$ $\theta_{n,m}=:\theta_{n+1}-\theta_{m+1}$ and $v_{n,m}:=v_n-v_m,$ then
$$
\left\lbrace
\begin{array}{l}
\partial_{t}\theta_{n,m}+v_n\cdot\nabla\theta_{n,m}+|\textnormal{D}|\theta_{n,m}=-v_{n,m}\cdot\nabla\theta_{m+1}\\
{\theta_{n,m}}_{|t=0}=0.
\end{array}
\right.
$$
Applying Theorem \ref{Thm3} to this equation gives
 \begin{equation}\label{Coupe1}
 \|\theta_{n,m}\|_{\widetilde L^\infty_t\dot{B}_{\infty,1}^0}
 \leq 
 Ce^{C\|\theta_n\|_{L^1_{t}\dot B^1_{\infty,1}}} \int_{0}^t\|v_{n,m}\cdot\nabla\theta_{m+1}(\tau)\|_{\dot{B}_{\infty,1}^0}d\tau.
 \end{equation}
Thanks to  Bony's decomposition \cite{bony},  the embedding $ \dot{B}^0_{\infty,1}\hookrightarrow L^\infty$ and the fact that ${\mathop{\rm div}}\,v_{n,m}=0$ 
\begin{equation}\label{Ac1}
\|v_{n,m}\cdot\nabla\theta_{m+1}\|_{\dot{B}_{\infty,1}^0}
\lesssim 
\|v_{n,m}\|_{\dot{B}_{\infty,1}^0}\|\theta_{m+1}\|_{\dot{B}_{\infty,1}^1}.
\end{equation}
Since Riesz transforms map continuously $\dot{B}_{\infty,1}^0$ into itself, then we get
\begin{equation}\label{Riesz}
\|v_{n,m}\cdot\nabla\theta_{m+1}\|_{\dot{B}_{\infty,1}^0}
\lesssim \|\theta_{n-1,m-1}\|_{\dot{B}_{\infty,1}^0}\|\theta_{m+1}\|_{\dot{B}_{\infty,1}^1}.
\end{equation}
Thus we infer
$$
 \|\theta_{n,m}\|_{\widetilde L^\infty_t\dot{B}_{\infty,1}^0}
 \leq 
 C \|\theta_{n-1,m-1}\|_{\widetilde L^\infty_t\dot{B}_{\infty,1}^0}
 e^{C\|\theta_n\|_{L^1_{t}\dot B^1_{\infty,1}}} \int_{0}^t\|\theta_{m+1}(\tau)\|_{\dot{B}_{\infty,1}^1}d\tau.
$$
According to the inequality (\ref{condition_petitesse}) one can choose $\varepsilon_0$ small such that
$$
 \|\theta_{n,m}\|_{\widetilde L^\infty_T\dot{B}_{\infty,1}^0}
 \leq 
 \eta\|\theta_{n-1,m-1}\|_{\widetilde L^\infty_T\dot{B}_{\infty,1}^0}
 $$
 with $\eta<1.$
Let us suppose that $n\geq m,$ then by induction  one finds
$$
 \|\theta_{n,m}\|_{\widetilde L^\infty_T\dot{B}_{\infty,1}^0}
 \lesssim
\eta^m\|\theta^0\|_{\dot{B}_{\infty,1}^0}.
$$
Thus $(\theta_n)_{n\in\NN}$ is a Cauchy sequence in $\widetilde L^\infty_T\dot{B}_{\infty,1}^0.$ Then there exists  $\theta 
\in\widetilde L^\infty_t\dot{B}_{\infty,1}^0$  such that $\theta_n$ converges strongly to $\theta$ in $\widetilde L^\infty_t\dot{B}_{\infty,1}^0.$
Moreover Fatou lemma and inequality (\ref{condition_petitesse}) imply that $\theta\in L^1_t\dot{B}_{\infty,1}^1.$
These informations allow us to pass to the limit into the equation.
\\

{\it Step 3:  Uniqueness}.\\
 Let us denote $X_{T}:=L^\infty_{T}\dot{B}_{\infty,1}^0\cap L^1_{T}\dot{B}_{\infty,1}^1$ and $\theta_{i}, i=1,2$ ( $v_{i}$ the corresponding velocity) be two solutions of the $\QG$ equation with the same  initial data  and  belonging to the space $X_{T}.$ We set $\theta_{1,2}=\theta_1-\theta_2$ and $v_{1,2}=v_1-v_2,$ then it is plain that
$$
\partial_{t}\theta_{1,2}+v^1\cdot\nabla\theta_{1,2}+|\textnormal{D}|\theta_{1,2}=-v_{1,2}\cdot\nabla\theta^{2},\,\,{\theta_{1,2}}_{|t=0}=0.
$$
Thanks to the inequalities (\ref{Coupe1}) and (\ref{Riesz}), we have  
$$
 \|\theta_{1,2}\|_{\widetilde L^\infty_t\dot{B}_{\infty,1}^0}
 \leq 
 Ce^{C\|\theta_1\|_{L^1_{t}\dot B^1_{\infty,1}}} \int_{0}^t\|\theta_{1,2}\|_{\widetilde L^\infty_{\tau}\dot B^0_{\infty,1}}
 \|\theta_2(\tau)\|_{\dot{B}_{\infty,1}^1}d\tau.
 $$
Thus  Gronwall's inequality gives the desired result.\\

{\it Step 4: smoothing effect.}

We will show the precise estimate: for all $\beta\in\RR_{+}$ we have
\begin{equation}\label{wassouf}
\| t^\beta\theta(t)\|_{\widetilde L^\infty_T\dot B^\beta_{\infty,1}}
\leq C_{\beta}
e^{C(\beta+1)\|\theta\|_{L^1_T\dot B^1_{\infty,1}}}\|\theta\|_{\widetilde L^\infty_T\dot B^0_{\infty,1}}.
\end{equation}
It is clear that
$$
\left\{
\begin{array}{rl}
&\partial_t(t^\beta\theta)+v\cdot\nabla(t^\beta\theta)+\vert \textnormal{D}\vert(t^\beta\theta)=\beta t^{\beta-1}\theta
\\
&(t^\beta\theta)_{|t=0}=0.
\end{array}
\right.
$$
We will proceed by induction and start the proof with the case $\beta\in\NN.$ 

For $\beta=1,$ we apply Theorem \ref{Thm3} with $\bar{r}=+\infty,$
$$
\| t\theta(t)\|_{\widetilde L^\infty_T\dot B^1_{\infty,1}}
\lesssim
e^{C\|\theta\|_{L^1_T\dot B^1_{\infty,1}}}\|\theta\|_{\widetilde L^\infty_T\dot B^0_{\infty,1}}.
$$
Assume (\ref{wassouf}) holds for degree $n;$ we will prove it for $n+1.$

Applying Theorem \ref{Thm3} to the equation of $t^{n+1}\theta$ we get
$$
\begin{aligned}
\| t^{n+1}\theta(t)\|_{\widetilde L^\infty_T\dot B^{n+1}_{\infty,1}}
&\leq C(n+1)
e^{C\|\theta\|_{L^1_T\dot B^1_{\infty,1}}}\|t^{n}\theta\|_{\widetilde L^\infty_T\dot B^n_{\infty,1}}
\\&
\leq C_{n}
e^{C(n+2)\|\theta\|_{L^1_T\dot B^1_{\infty,1}}}\|\theta\|_{\widetilde L^\infty_T\dot B^0_{\infty,1}}.
\end{aligned}
$$
For $\beta\in\RR_+,$ we have $[\beta]\leq\beta<[\beta]+1$ and by interpolation, one has
$$
\begin{aligned}
\|t^\beta\theta\|_{\widetilde L^\infty_T(\dot B^\beta_{\infty,1})}
&\lesssim
\Big\|t^{[\beta]}\theta\Big\|_{\widetilde L^\infty_T(\dot B^{[\beta]}_{\infty,1})}^{1+[\beta]-\beta}
\Big\|t^{[\beta]+1}\theta\Big\|_{\widetilde L^\infty_T(\dot B^{[\beta]+1}_{\infty,1})}^{\beta-[\beta]}.
\end{aligned}
$$
This completes the proof.
\end{proof}

\subsection{Blowup Criteria}
The main result of this paragraph is:
\begin{prop}\label{critere_explosion}
Let $T^*$ be the maximum local existence time of $\theta$ in $\widetilde L^\infty_T\dot B^0_{\infty,1}\cap 
L^1_T\dot B^1_{\infty,1}.$ There exists an absolute constant $\varepsilon_0>0$ such that \mbox{if $T^*<\infty,$} then
$$
\displaystyle\liminf_{t \to T^*}(T^*-t)\|\nabla\theta(t)\|_{L^\infty}\geq\varepsilon_0.
$$ 
\end{prop}
\begin{proof}
From  local existence theory and especially (\ref{temps_existence}) we see that \mbox{if  $T^*<\infty,$} then necessary
$$
\displaystyle\liminf_{t \to T^*}
\sum_{q\in\ZZ}(1-e^{-c(T^*-t)2^q})^{\frac{1}{2}}\|\theta_q(t)\|_{L^\infty}
\geq
\varepsilon_0,
$$
otherwise we can continue the solution over $T^{*}.$
It follows that
$$
\displaystyle\liminf_{t \to T^*}\sum_{q\in\ZZ}(1-e^{-c(T^*-t)2^q})^{\frac{1}{2}}\displaystyle\sup_{t \leq T^*}\|\theta_q(t)\|_{L^\infty}
\geq
\varepsilon_0,
$$
Consequently we obtain from Lebesgue theorem
$$
\|\theta\|_{\widetilde L^\infty_{T^*}(\dot B^0_{\infty,1})}=\infty.
$$
Using Bernstein inequality and the fact that $\|\theta_q\|_{L^\infty}\lesssim\|\theta^0\|_{L^\infty},$
we have
$$
\begin{aligned}
\varepsilon_0
\leq
\displaystyle\liminf_{t \to T^*}&\Big\{\sum_{q\leq N}(1-e^{-c(T^*-t)2^{q}})^{\frac{1}{2}}\|\theta_q(t)\|_{L^\infty}
\\&
+\sum_{q\geq N}(1-e^{-c(T^*-t)2^q})^{\frac{1}{2}}
\|\theta_q(t)\|_{L^\infty}\Big\}
\\&
\lesssim
\displaystyle\liminf_{t \to T^*}\Big\{(T^*-t)^{\frac{1}{2}}\|\theta^0\|_{L^\infty} \sum_{q\leq N}2^{q/2}
+ \|\nabla\theta(t)\|_{L^\infty}\sum_{q\geq N}2^{-q}\Big\}
\\&
\lesssim
\displaystyle\liminf_{t \to T^*}\Big\{(T^*-t)\|\theta^0\|_{L^\infty} 2^N
+ \|\nabla\theta(t)\|_{L^\infty}2^{-N}\Big\}.
\end{aligned}
$$
Choosing judiciously $N$ we obtain the desired result.
\end{proof}
\subsection{Global existence.}
We will use the idea of \cite{volberg}.
 Let $T^*$ be the maximal time existence of the solution in the space $\widetilde L^\infty_{\textnormal{loc}}([0,T^*[, \dot{B}_{\infty,1}^0)
\cap L^1_{\textnormal{loc}}([0,T^*[, \dot{B}_{\infty,1}^1).$ From the local existence, there exists $T_{0}>0$ such that 
 $$
\forall t\in[0,T_{0}],\,t\|\nabla\theta(t)\|_{L^\infty}
\leq 
C\|\theta^0\|_{\dot{B}_{\infty,1}^0}.
$$
Let $\lambda$ be a real positive number that will be fixed later and  $T_{1}\in]0,T_{0}[.$ We define
the set
$$
I:=\Big\{T\in[T_{1},T^*[; \forall t\in [T_{1},T], \forall x\neq y\in\RR^2, |\theta(t,x)-\theta(t,y)|<\omega_{\lambda}(|x-y|)\Big\},
$$
where
$$
\omega:\hspace{0,2cm} \RR_+ \longrightarrow \RR_+,
$$
is strictly increasing, concave, $\omega(0)=0,$ $\omega'(0)<+\infty,$
$\displaystyle\lim_{\xi\longrightarrow 0^{+}}\omega''(\xi)=-\infty$
and
$$
\omega_{\lambda}(|x-y|)=\omega(\lambda|x-y|).
$$
The function $\omega$ is a modulus of continuity chosen as in \cite{volberg}. We shall first check that $I$ is nonempty. It suffices for this purpose to  prove that $T_{1}$ belongs to $I$ under suitable conditions over $\lambda.$ Let $C_{0}$ be a large positive number such that 
 \begin{equation}\label{R1}
 \omega(C_{0})> 2\|\theta^0\|_{L^\infty}.
 \end{equation}
 Since $\omega$ is strictly non-decreasing function then we get from maximum principle
 $$
\forall x,y; \lambda|x-y|\geq C_{0}\Rightarrow 
|\theta(T_{1},x)-\theta(T_{1},y)|\leq2\|\theta^0\|_{L^\infty}
<\omega_{\lambda}(|x-y|).
$$ 
On the other hand we have  from Mean Value Theorem 
$$
|\theta(T_{1},x)-\theta(T_{1},y)|
\leq 
|x-y|\|\nabla \theta(T_{1})\|_{L^\infty}. 
$$
Let $0<\delta_{0}<C_{0}.$ Then using the concavity of $\omega$ one obtains
$$
\lambda|x-y|\leq\delta_{0}\Rightarrow\omega_{\lambda}(|x-y|)\geq\frac{\omega(\delta_{0})}{\delta_{0} }\lambda|x-y|.
$$
If we choose $\lambda$ so that 
$$
\lambda>\frac{\delta_{0}}{\omega(\delta_{0})}
\|\nabla\theta(T_{1})\|_{L^\infty},$$ then we get
$$
0<\lambda|x-y|\leq\delta_{0}\Rightarrow 
|\theta(T_{1},x)-\theta(T_{1},y)|<\omega_{\lambda}(|x-y|).
$$
Let us now move to the case $\delta_{0}\leq \lambda|x-y|
\leq C_{0}.$ By an obvious computation we find
\begin{eqnarray*}
|\theta(T_{1},x)-\theta(T_{1},y)|&\leq& \frac{C_{0}}{\lambda}\|\nabla\theta(T_{1})\|_{L^\infty}\quad\hbox{and}\\
\omega(\delta_{0})&\leq&\omega(\lambda|x-y|).  
\end{eqnarray*}
Choosing $\lambda$ such that
$$
\lambda>\frac{C_0}{\omega(\delta_{0})}
\|\nabla\theta(T_{1})\|_{L^\infty}.
$$
Then we obtain
$$
 \delta_{0}\leq \lambda|x-y|\leq C_{0}\Rightarrow 
 |\theta(T_{1},x)-\theta(T_{1},y)|<\omega_{\lambda}(|x-y|).
 $$
All the preceding conditions over $\lambda$ can be obtained if we take
\begin{equation}\label{abas}
\lambda=\frac{\omega^{-1}(3\|\theta^0\|_{L^\infty})}
{2\|\theta^0\|_{L^\infty}}\|\nabla\theta(T_{1})\|_{L^\infty}.
\end{equation}
From the construction, the set $I$ is an interval of the form $[T_{1},T_{*}).$ We have three possibilities.
The first one is $T_{*}=T^*$ and in this case we have necessary $T^*=+\infty$ because the Lipshitz norm of $\theta $ does not blow up. The second one is $T_{*}\in I$   and we will show that is not possible. Indeed,
let $C_{0}$ as (\ref{R1}) then for all $t\in[T_{1},T^*)$ 
$$
\lambda|x-y|\geq C_{0}\Rightarrow|\theta(t,x)-\theta(t,y)|<\omega_{\lambda}(|x-y|)
$$
Since $\nabla\theta(t)$ belongs to $C(]0,T^*);\dot{B}_{\infty,1}^0),$ then for $\epsilon>0$ there exist $\eta_{0}, R>0$ such that
$$
\forall t\in[T_{*},T_{*}+\eta_{0}]\Rightarrow \|\nabla\theta(t)\|_{L^\infty}\leq \|\nabla\theta(T_{*})\|_{L^\infty}+\epsilon/2\quad\hbox{and}
$$
$$
\|\nabla\theta(T_{*})\|_{L^\infty(B_{(0,R)}^c)}\leq \epsilon/2,
$$
where $B_{(0,R)}$ is the ball of radius $R$ and with center the origin.

Hence for $\lambda|x-y|\leq C_{0}$ and $x$ or $y\in B_{(0,R+\frac{C_{0}}{\lambda})}^c$ we have for $t\in[T_{*},T_{*}+\eta_{0}]$
 \begin{eqnarray*}
 |\theta(t,x)-\theta(t,y)|&\leq&|x-y|\|\nabla\theta (t)\|_{L^\infty(B_{(0,R)}^c)}\\
 &\leq& \epsilon|x-y|.
 \end{eqnarray*}
 On  the other hand we have from the concavity of $\omega$
 $$
 \lambda|x-y|\leq C_{0}\Rightarrow\frac{\omega(C_{0})}{C_{0}}\lambda|x-y|\leq\omega_{\lambda}(|x-y|)
 $$
Thus if we take $\epsilon$ sufficiently small such that
$$
\epsilon<\frac{\omega(C_{0})}{C_{0}}\lambda,$$
then we find that 
$$
\lambda|x-y|\leq C_{0};x\,\hbox{or}\,y\in B_{(0,R+\frac{C_{0}}{\lambda})}^c\Rightarrow|\theta(t,x)-\theta(t,y)|<\omega_{\lambda}(|x-y|).
$$
It remains to study the case where $x,y\in B_{(0,R+\frac{C_{0}}{\lambda})}.$ Since $\|\nabla^2\theta(T_{*})\|_{L^\infty}$is finite (see Proposition \ref{N67}) then we get for each $x\in\RR^2$
$$
|\nabla\theta(T_{*},x)|<\lambda\omega'(0).
$$
For the proof see \cite{Yu}, page $4.$
From the continuity of $x\longmapsto|\nabla\theta(T_{*},x)|$ we obtain
$$
\|\nabla\theta(T_{*})\|_{L^\infty(B_{(0,R+\frac{C_{0}}{\lambda})})}<\lambda\omega'(0)
$$
Let $\delta_{0}<<1$ then using the continuity in time of the quantity $\|\nabla\theta(t)\|_{L^\infty}$ one can find $\eta_{1}>0$ such that $\forall t\in[T_{*}, T_{*}+\eta_{1}]$
$$
\|\nabla\theta(t)\|_{L^\infty(B_{(0,R+\frac{C_{0}}{\lambda})})}< \lambda\frac{\omega(\delta_{0})}{\delta_{0}}\cdot$$
For $\lambda|x-y|\leq\delta_{0}$ and $x\neq y$ belonging together to $ B_{(0,R+\frac{C_{0}}{\lambda})}$ we  have
\begin{eqnarray*}
|\theta(t,x)-\theta(t,y)|&\leq &|x-y|\|\nabla\theta(t)\|_{L^\infty(B_{(0,R+\frac{C_{0}}{\lambda})})}\\
&<& \lambda|x-y|\frac{\omega(\delta_{0})}{\delta_{0}}\leq \omega_{\lambda}(|x-y|).
\end{eqnarray*}
Now for the other case we have
$$
\forall x,y\in B_{(0,R+\frac{C_{0}}{\lambda})}, \delta_{0}\leq \lambda|x-y|; |\theta(T_{*},x)-\theta(T_{*},y)|<\omega_{\lambda}(|x-y|), 
$$
then we get from a standard compact argument the existence of $\eta_{2}>0$ such that for all $t\in[T_{*},T_{*}+\eta_{2}]$
$$
\forall x,y\in B_{(0,R+\frac{C_{0}}{\lambda})}, \delta_{0}\leq \lambda|x-y|; |\theta(t,x)-\theta(t,y)|<\omega_{\lambda}(|x-y|).$$
Taking $\eta=\min\{\eta_{0},\eta_{1},\eta_{2}\},$ we obtain that $T_{*}+\eta\in I$ which contradicts the fact that $T_{*}$ is maximal.

The last case that we have to treat is that $T_{*}$ does not belong to $I.$ Thus we have by time continuity of $\theta$ the existence of $x\neq y$ such that
$$
\theta(T_{*},x)-\theta(T_{*},y)=\omega_{\lambda}(\xi), \hbox{with}\quad\xi=|x-y|.
$$
We will show that this scenario can not occur and more precisely:  
$$
f'(T_{*})<0\,\quad\hbox{where}\quad f(t)=\theta(t,x)-\theta(t,y).
$$
This is impossible since $f(t)\leq f(T_{*}), \forall t\in[0,T_{*}].$
The proof is the same as \cite{volberg} and for the convenience of the reader we will outline the proof.
From the regularity of the solution we see that the $\QG$ equation can be defined in the classical manner and
$$
f'(T_{*})=(u\cdot\nabla \theta)(T_{*},x)-(u\cdot\nabla \theta)(T_{*},y)+|\hbox{D}|\theta(T_{*},x)-|\hbox{D}|\theta(T_{*},y).
$$
From \cite{volberg} we have
$$
(u\cdot\nabla \theta)(T_{*},x)-(u\cdot\nabla \theta)(T_{*},y)\leq \Omega_{\lambda}(\xi)\omega_{\lambda}'(\xi),$$
where
$$
\Omega_{\lambda}(\xi)=C\Big(\int_{0}^\xi\frac{\omega_{\lambda}(\eta)}{\eta}d\eta+\xi\int_{\xi}^\infty\frac{\omega_{\lambda}(\eta)}{\eta^2}d\eta\Big)=\Omega(\lambda\xi).
$$
Again from \cite{volberg}
\begin{eqnarray*}
|\hbox{D}|\theta(T_{*},x)-|\hbox{D}|\theta(T_{*},y)&\leq&\frac{1}{\pi }\int_{0}^{\frac{\xi}{2}}\frac{\omega_{\lambda}(\xi+2\eta)+\omega_{\lambda}(\xi-2\eta)-2\omega_{\lambda}(\xi)}{\eta^2}d\eta\\
&+&\frac{1}{\pi }\int_{\frac{\xi}{2}}^\infty\frac{\omega_{\lambda}(2\eta+\xi)-\omega_{\lambda}(2\eta-\xi)-2\omega_{\lambda}(\xi)}{\eta^2}d\eta\\
&\leq&\lambda\, {\mathcal I}(\lambda\xi),
\end{eqnarray*}
 where
\begin{eqnarray*}
{\mathcal{I}}(\xi)&=&\frac{1}{\pi }\int_{0}^{\frac{\xi}{2}}\frac{\omega(\xi+2\eta)+\omega(\xi-2\eta)-2\omega(\xi)}{\eta^2}d\eta\\
&+&\frac{1}{\pi }\int_{\frac{\xi}{2}}^\infty\frac{\omega(2\eta+\xi)-\omega(2\eta-\xi)-2\omega(\xi)}{\eta^2}d\eta.
\end{eqnarray*}
Thus we get
$$
f'(T_{*})=\lambda(\Omega \omega'+{\mathcal{I}})(\lambda\xi).
$$
Now, we choose the same function as \cite{volberg} (see page 5)
$$
\omega(\xi)=\xi-\xi^{\frac{3}{2}}, \hbox{if}\quad \quad\xi\in[0,\delta]
$$
and
$$
\omega'(\xi)=\frac{\gamma}{\xi(4+\log(\xi/\delta))}, \hbox{if}\quad\quad \xi>\delta,
$$
where $\delta$ and $\gamma$ are small numbers and satisfy $0<\gamma<\delta.$
It is shown in \cite{volberg} that
$$
\Omega(\xi) \omega'(\xi)+{\mathcal{I}}(\xi)<0, \forall\xi\neq 0.
$$
This yields to
 $f'(T_{*})<0.
 $
 
Finally we have $T^*=+\infty$ and 
$$
\forall t\in[T_{1},+\infty), \,\,\|\nabla\theta(t)\|_{L^\infty}\leq \lambda.
$$
The value of $\lambda $ is given by (\ref{abas}).

\end{document}